\documentclass[11pt]{paper}
\usepackage{amssymb,euscript,graphicx,pifont,cite,amsmath,epsfig}
\usepackage{cite}
\usepackage{pifont}
\usepackage{times}
\usepackage{amssymb,epsfig,floatflt}
\usepackage{subfigure}[1995/03/06]
\usepackage[section]{placeins}
\usepackage{epsfig,psfrag}
\usepackage{graphics}
\usepackage{color}

\input amssymb.sty

\def\dspace{\multiply\normalbaselineskip 150

          \divide\normalbaselineskip 100 \normalbaselines

\csname @@normalbaselineskip\endcsname\normalbaselineskip}

\def\rr{\rm I\!R}

\newtheorem{thm}{Theorem}

\newtheorem{lem}{Lemma}

\newtheorem{defn}{Definition}

\newtheorem{cor}{Corollary}

\newtheorem{rem}{Remark}

\usepackage{cite}

\begin{document}
\thispagestyle{empty}

\title{\LARGE \bf Guaranteed Transient Performance with $\mathcal{L}_1$ Adaptive Controller
for Systems with Unknown Time-varying  Parameters: Part I \thanks{Research is supported by AFOSR under Contract
No. FA9550-05-1-0157.}}
\author{Chengyu Cao and Naira Hovakimyan\thanks{The authors are with
  Aerospace \& Ocean Engineering, Virginia Polytechnic
 Institute \& State University, Blacksburg, VA 24061-0203, e-mail:
  chengyu, nhovakim@vt.edu}}
\date{}
\maketitle
\begin{abstract}
This paper presents a novel adaptive control methodology for
uncertain systems with time-varying unknown parameters and
time-varying bounded disturbance. The adaptive controller
ensures uniformly bounded transient and asymptotic tracking for
system's both signals, input and output, simultaneously. The
performance bounds can be systematically improved by increasing
the adaptation gain. Simulations of a robotic arm with
time-varying friction verify the theoretical findings.
\end{abstract}

\section{Introduction}

This paper presents an adaptive control methodology for
controlling systems with unknown time-varying parameters, which
are not required to have slow variation. The methodology ensures
uniformly bounded transient response for system's both signals,
input and output, simultaneously, in addition to asymptotic
tracking. The main advantage of this new architecture, as compared
to the existing results in the literature,  is that it ensures
uniform transient tracking for system's input signal in addition
to its output. The $\mathcal L_\infty$ norm bounds for the error
signals between the closed-loop adaptive system and  the
closed-loop reference LTI system can be systematically reduced by
increasing the adaptation gain.

Adaptive algorithms  achieving arbitrarily improved transient
performance in case of constant unknown parameters are given in
\cite{DattaHo94, BartFerSto99, Sun93, MilTac91acw, Costa99,
YdsTac92tpa, KrsScl93tpi, OrtTac93mna, ZanCdc90tbf, Datta94,
arteaga2002, narendra94}, and for unknown time-varying parameters
have been given in \cite{marinotomei99, marinotomei99tac}. While
the results in \cite{marinotomei99, marinotomei99tac} improved
upon \cite{Tsakalis87, middleton88, zang96}, by extending the
class of systems beyond the slow time-variation of the unknown
parameters and guaranteeing performance improvement to arbitrary
degree, they still did not provide means for regulating the
 performance of the control signal during the transient.

A common tendency observed in a variety of applications using
adaptive control is that increasing the adaptation gain leads to
improved transient tracking of the system output, but the  control
signal experiences  high-frequency oscillations. In \cite{kkk}, a
bound is derived to confirm the first part of this statement
assuming appropriate trajectory initialization.  The
high-frequency oscillations in the control signal consequently
limit the rate of adaptation. If one considers the simplest
adaptive scheme for a scalar linear system with constant
disturbance, which can be solved by a PI controller, then it is
straightforward to verify that increasing the adaptation gain
leads to reduced phase margin for the resulting closed-loop linear
system, \cite{chengyu_gnc06}. This observation explains to some
extent the oscillations inherent to the control signal in the
presence of high adaptation gain.

In recent papers \cite{CaoHovakimyan,CaoHovakimyan1}, we have
developed a novel $\mathcal{L}_1$ adaptive control architecture
that permits fast adaptation and yields guaranteed transient
response for system's both signals, input and output,
simultaneously, in addition to asymptotic tracking. The  main
feature of it is the ability of fast adaptation with guaranteed
low-frequency control signal. The ability of fast adaptation
ensures the desired transient performance for system's both
signals, input and output, simultaneously, while the low-pass
filter in the feedback loop attenuates the high-frequency
components in the control signal.  In this paper we expand the
class of systems to have time-varying unknown parameters of
arbitrary rate of variation, and we correspondingly modify the
architecture from \cite{CaoHovakimyan,CaoHovakimyan1} to ensure
the desired transient performance for system's both signals. We
prove that by increasing the adaptation gain one can achieve
arbitrary close transient and asymptotic tracking for system's
both signals, input and output, simultaneously. In Part II of this
paper\cite{CDC06_chengyu_big_L1}, we prove that increasing the
adaptation gain will not hurt the time-delay margin of the
closed-loop system with the $\mathcal{L}_1$ adaptive control
architecture, as opposed to the conventional adaptive schemes
observed in \cite{chengyu_gnc06}.

The paper is organized as follows. Section \ref{sec:preliminary}
states some preliminary definitions, and Section \ref{sec:PF}
gives the problem formulation. In Section \ref{sec:Hinf}, the
novel $\mathcal{L}_1$ adaptive control architecture  is presented.
Stability and uniform transient tracking bounds of the
$\mathcal{L}_1$ adaptive controller are presented in Section
\ref{sec:convergence}. In section \ref{sec:simu}, simulation
results are presented, while Section \ref{sec:con} concludes the
paper.

\section{Preliminaries}\label{sec:preliminary}

In this Section, we recall some basic definitions and facts from
linear systems theory, \cite{IoaBook03,KhaBook02,ZhoBook98}.

%
%
%
%
%
%
%
%
%
%

\begin{defn}\label{defn1}
For a signal $\xi(t), ~t\geq 0, ~\xi\in {\rr}^n$, its truncated
${\mathcal L}_{\infty}$ norm and ${\mathcal L}_{\infty}$ norm are
defined as
\begin{eqnarray}
 \Vert \xi_t \Vert_{{\mathcal L}_\infty} & = &
\max_{i=1,..,n} \Big(\sup_{0\leq \tau \leq t}
|\xi_i(\tau)|\Big)\,,\nonumber\\
\Vert \xi \Vert_{{\mathcal L}_\infty} & = &  \max_{i=1,..,n}
\Big(\sup_{\tau\geq 0} |\xi_i(\tau)|\Big)\,,\nonumber
\end{eqnarray}
 where $\xi_i$ is the
$i^{th}$ component of $\xi$.
\end{defn}
\begin{defn}
 The $\mathcal{L}_1$ gain of a stable proper single--input single--output system $H(s)$
 is defined to be
$|| H(s)||_{\mathcal{L}_1} = \int_{0}^{\infty} |h(t)| d t, $ where
$h(t)$ is the impulse response of $H(s)$, computed via the inverse
Laplace transform $ h(t)=\frac{1}{2 \pi i}\int_{\alpha-i
\infty}^{\alpha+i \infty} H(s) e^{st} ds, t\geq 0, $ in which the
integration is done along the vertical line $x=\alpha>0$ in the
complex plane.
\end{defn}

{\em Proposition:}
 A continuous time LTI
system (proper) with impulse response $h(t)$ is stable if and only
if
$
\int_{0}^{\infty} |h(\tau)| d\tau < \infty.
$
A proof can be found in \cite{IoaBook03} (page 81, Theorem 3.3.2).

\begin{defn}
For a stable proper $m$ input $n$ output system $H(s)$ its
$\mathcal{L}_1$ gain is defined as
\begin{equation}\label{L1def}
\Vert H(s) \Vert_{\mathcal{L}_1} = \max_{i=1,\cdots,n} \left(
\sum_{j=1}^{m}\Vert H_{ij}(s)\Vert_{\mathcal{L}_1} \right)\,,
\end{equation}
where $H_{ij}(s)$ is the $i^{th}$ row $j^{th}$ column element of
$H(s)$.
\end{defn}

The next lemma extends the results of Example 5.2
(\cite{KhaBook02}, page 199)  to general multiple input multiple
output systems.

\begin{lem}\label{lem:L1} For a stable proper multi-input multi-output (MIMO)  system $H(s)$ with input $r(t) \in {\rr}^m$ and output
 $x(t)\in {\rr}^n$, we have
\begin{equation}
\Vert x_t\Vert_{{\mathcal L}_\infty} \leq \Vert H\Vert_{{\mathcal
L}_1} \Vert r_t\Vert_{{\mathcal L}_\infty},\quad \forall~
t>0.\nonumber
\end{equation}
\end{lem}
%
%
%

\begin{cor}\label{lem:L1_ext} For a stable proper MIMO  system  $H(s)$, if the input   $r(t) \in {\rr}^m$ is bounded, then the  output
 $x(t)\in {\rr}^n$ is also bounded as
$ 
\Vert x\Vert_{{\mathcal L}_\infty } \leq \Vert
H(s)\Vert_{\mathcal{L}_1} \Vert r\Vert_{{\mathcal L}_\infty }.
$ 
\end{cor}

\begin{lem}\label{lem:L1cas} For a cascaded system $H(s)=H_2(s) H_1(s)$,
where $H_1(s)$ is a stable proper system with $m$ inputs and $l$
outputs
 and $H_2(s)$ is a stable   proper system with $l$ inputs and  $n$
outputs, we have
$ 
\Vert H(s) \Vert_{\mathcal{L}_1} \leq  \Vert H_2(s)
\Vert_{\mathcal{L}_1} \Vert H_1(s) \Vert_{\mathcal{L}_1}\,.
$
\end{lem}

\begin{figure}[!t]
\begin{center}
\includegraphics[width=2.5in,height=1.4in]{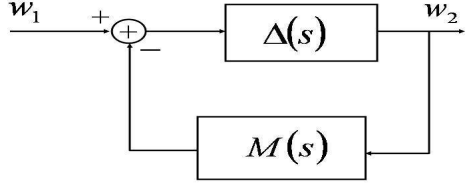}
\caption{Interconnected systems}\label{fig:sgblock}
\end{center}
\end{figure}

 Consider an interconnected
LTI system in Fig. \ref{fig:sgblock}, where $w_1\in {\rr}^{n_1}$,
$w_2\in {\rr}^{n_2}$, $M(s)$ is a stable  proper system with $n_2$
inputs and $n_1$ outputs, and $\Delta(s)$ is a stable proper
system with $n_1$ inputs and $n_2$ outputs.
\begin{thm}\label{thm:sg}
({\bf $\mathcal{L}_1$ Small Gain Theorem})  The interconnected
system in Fig. \ref{fig:sgblock} is stable if
$ 
\Vert M(s) \Vert_{\mathcal{L}_1} \Vert \Delta(s)
\Vert_{\mathcal{L}_1} < 1.
$ 
\end{thm}
The proof follows from Theorem 5.6  (\cite{KhaBook02}, p. 218),
written for $\mathcal{L}_1$ gain.

Consider a linear time invariant system:
\begin{equation}\label{LTI_def}
\dot{x}(t) = A x(t) +b u(t)\,,
\end{equation}
where $x\in {\rr}^n$, $u\in {\rr}$, $b\in {\rr}^{n}$, $A\in
{\rr}^{n \times n}$ is Hurwitz, and assume that the transfer
function $(s I -A)^{-1} b $ is strictly proper and stable. Notice
that it can be expressed as:
\begin{equation}\label{LTI_trans}
(s I -A)^{-1} b = \frac{n(s)}{d(s)}\,,
\end{equation}
where $d(s)={\rm{det}} (s I- A)$ is a $n^{th}$ order stable
polynomial, and $n(s)$ is a $n\times 1$ vector with its $i^{th}$
element being a polynomial function:
\begin{equation}\label{nijdef}
n_i(s) = \sum_{j=1}^{n} n_{ij} s^{j-1} \,.
\end{equation}

\begin{lem}\label{lem:pre1}
If $(A\in {\rr}^{n\times n}, b\in {\rr}^{n})$ is controllable, the
matrix $N$ with its $i^{th}$ row $j^{th}$ column entry $n_{ij}$ is
full rank.
\end{lem}

\begin{lem}\label{lem:pre2}
If $(A, b)$ is controllable and $(s I -A)^{-1} b $ is strictly
proper and stable, there exists $c\in {\rr}^n$ such that the
transfer function $c^{\top}(s I -A)^{-1} b$ is minimum phase
 with relative degree one, i.e. all its zeros are located in
the left half plane, and its denominator is one order larger than
its numerator.
\end{lem}

\section{Problem Formulation}\label{sec:PF}
Consider the following system dynamics:
\begin{eqnarray}
\dot{x}(t)  & = &  A_m x(t)+b \left( \omega u(t)+\theta^{\top}(t)
x(t)+ \sigma(t) \right)\,,\nonumber\\
y(t) & = & c^{\top} x(t),\quad x(0)=x_0\,,\label{problemnew}
\end{eqnarray}
 where $x\in
{\rr}^n$ is the system state vector (measurable), $u\in {\rr}$ is
the control signal, $y\in {\rr}$ is the regulated output, $b,c\in
{\rr}^n$ are known constant vectors, $A_m$ is a known $n \times n$
matrix,  $\omega\in \rr$ is  an unknown constant with known sign,
$\theta(t)\in \rr^n$ is a vector of time-varying unknown
parameters, while $\sigma(t)\in \rr$ is a time-varying
disturbance. Without loss of generality, we assume that
\begin{equation}\label{Thetadef}
\omega\in \Omega=[\omega_{l},\;\omega_{u}]\,,\theta(t)\in
\Theta,\; |\sigma(t)|\leq \Delta\,, \quad t\geq 0\,,
\end{equation}
where $\omega_{u}>\omega_{l}>0$ are given bounds, $\Theta$ is
known compact set and $\Delta\in \rr^{+}$ is a known
(conservative) $\mathcal{L}_{\infty}$ bound of $\sigma(t)$.

The control objective is to design a full-state feedback adaptive
controller to ensure that $y(t)$ tracks a given bounded reference
signal $r(t)$ {\em both in transient and steady state}, while all
other error signals remain bounded.

We further assume that $\theta(t)$ and $\sigma(t)$ are
continuously differentiable and their derivatives are uniformly
bounded:
\begin{equation}
\Vert \dot{\theta}(t) \Vert_{2}  \leq  d_{\theta}<\infty,\quad
|\dot{\sigma}(t)|  \leq  d_{\sigma}<\infty,\, \quad \forall~ t\geq
0\,
 ,\label{derivativeB}
\end{equation}
where $\Vert \cdot \Vert_2$ denotes the $2$-norm, while the
numbers $d_{\theta}, d_{\sigma}$ can be arbitrarily large.

\section{ $\mathcal{L}_1$ Adaptive Controller}\label{sec:Hinf}

In this section, we develop a novel adaptive control architecture
for the system in (\ref{problemnew}) that permits complete
transient characterization for both $u(t)$ and $x(t)$.
The elements of $\mathcal{L}_1$ adaptive controller are introduced
next:

{\bf Companion Model:} We consider the following companion model:
\begin{eqnarray}
\dot{\hat{x}}(t) & = & A_m \hat{x}(t)+b \left( \hat{\omega}(t)
u(t) + \hat{\theta}^{\top}(t) x(t)+\hat{\sigma}(t) \right)\,,
\nonumber\\
\hat{y}(t) & = & c^{\top} \hat{x}(t)\,,\quad \hat
x(0)=x_0\,,\label{L1_companionmodal}
\end{eqnarray}
which has the same  structure as the system in (\ref{problemnew}).
The only difference is that the  unknown parameters $\omega,
\theta(t), \sigma(t)$ are replaced by their adaptive estimates
$\hat{\omega}(t), \hat{\theta}(t), \hat{\sigma}(t)$ that are
governed by the following adaptation laws.

{\bf Adaptive Laws:} Adaptive estimates are given by:
\begin{eqnarray}
\dot{\hat{\theta}}(t) &  = & \Gamma_{\theta}{\rm Proj}(-x(t)
\tilde{x}^{\top}(t) P b , \hat{\theta}(t)), ~
{\hat{\theta}}(0)=\hat \theta_0  \label{adaptivelaw_L11}\\
\dot{\hat{\sigma}}(t) &  = & \Gamma_{\sigma}{\rm Proj}(
-\tilde{x}^{\top}(t) P b ,  \hat{\sigma}(t)), ~~~~~
{\hat{\sigma}}(0)=\hat \sigma_0  \label{adaptivelaw_L12}\\
\dot{\hat{\omega}}(t) &  = & \Gamma_{\omega}{\rm Proj}(
-\tilde{x}^{\top}(t) P b u(t), \hat{\omega}(t)),
{\hat{\omega}}(0)=\hat \omega_0 \label{adaptivelaw_L13}
\end{eqnarray}
where $\tilde{x}(t)=\hat x(t)-x(t)$ is the error signal between
the state of the system and the companion model, $\Gamma_{\theta}
=\Gamma_c I_{n\times n}\in {\rr}^{n\times n}$,
$\Gamma_{\sigma}=\Gamma_{\omega}=\Gamma_c$ are adaptation gains
with $\Gamma_c\in \rr^{+}$, and $P$ is the solution of the
algebraic equation $ A_m^{\top} P+P A_m =- Q$, $Q>0$.

{\bf Control Law:} The control signal is generated through gain
feedback of the following  system:
\begin{eqnarray}
\chi(s) & = & D(s) r_u(s) \,,\nonumber\\
u(s) & = & - k \chi(s)\,,\label{controllaw}
\end{eqnarray}
where $r_u(s)$ is the Laplace transformation of
$r_u(t)=\hat{\omega}(t) u(t)+\bar{r}(t)$,
\begin{equation}\label{barrt_def}
\bar{r}(t)=\hat{\theta}^{\top}(t) x(t)+\hat{\sigma}(t)-k_g r(t),
\end{equation}
\begin{equation}\label{kg_def}
 k_g =-\frac{1}{ c^{\top} A_m^{-1} b}\,,
\end{equation}
$k\in \rr^{+}$ is a
feedback gain, while $D(s)$ is any transfer function that leads to
strictly proper stable
\begin{equation}\label{Csdef}
C(s) = \frac{\omega k D(s)}{1+\omega k D(s)}\,
\end{equation}
 with low-pass gain $C(0)=1$.
 One simple choice is
\begin{equation}\label{1stDs}
D(s) =\frac{1}{s}\,,
\end{equation}
which yields a   first order strictly proper $C(s)$ in the
following form:
\begin{equation}\label{1stCs}
C(s) = \frac{\omega k }{s+\omega k }\,.
\end{equation}
Further, let
\begin{equation}\label{thetabardef} L =
\max_{\theta(t)\in \Theta} \sum_{i=1}^{n} |\theta_i(t)|\,,
\end{equation}
where $\theta_i(t)$ is the $i^{th}$ element of $\theta(t)$,
$\Theta$ is the compact set defined in (\ref{Thetadef}). We now
state the $\mathcal{L}_1$ performance requirement that ensures
 stability of the entire system and desired transient performance, as discussed later in Section \ref{sec:convergence}.

{\bf $\mathcal{L}_1$-gain stability requirement:} Design $D(s)$ to
ensure that
\begin{equation}\label{condition3}
\Vert G(s) \Vert_{\mathcal{L}_1} L <1\,,
\end{equation}
where $G(s)=(sI-A_m)^{-1}b (1-C(s))$.

\begin{figure}[!h]
\begin{center}
\includegraphics[width=3.4in,height=2.7in]{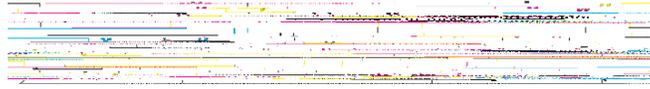}
\caption{Closed-loop system with $\mathcal{L}_1$ adaptive
controller}\label{fig:L1block}
\end{center}
\end{figure}

The complete $\mathcal{L}_1$ adaptive controller consists of
(\ref{L1_companionmodal}),
(\ref{adaptivelaw_L11})-(\ref{adaptivelaw_L13}) and
(\ref{controllaw}) subject to $\mathcal{L}_1$-gain stability
requirement in (\ref{condition3}).  The closed-loop system is
illustrated in Fig. \ref{fig:L1block}.

In case of constant $\theta(t)$,  the stability requirement of the
$\mathcal{L}_1$ adaptive controller can be simplified. For the
specific choice of  $D(s)$ and $C(s)$ in (\ref{1stDs}) and
(\ref{1stCs}), the stability requirement of $\mathcal{L}_1$
adaptive controller is reduced to
\begin{equation}
A_{g}=\left[
\begin{array}{cc}
A_m+b \theta^{\top} & b \omega\\
-k \theta^{\top} & -k \omega
\end{array}
\right]\, \label{condition4}
\end{equation}
being Hurwitz for all $\theta\in \Theta$, $\omega\in \Omega$.

\section{Analysis of $\mathcal{L}_1$ Adaptive
Controller}\label{sec:convergence}

\subsection{Closed-loop Reference System}
We now consider the following closed-loop LTI reference system
with its control signal  and system response being defined as
follows:
\begin{eqnarray}
\dot{x}_{ref}(t) & = & A_m x_{ref}(t)+\nonumber\\
& & b \left( \omega u_{ref}(t)+\theta^{\top}(t) x_{ref}(t)+\sigma(t)\right), \label{refx}\\
u_{ref}(s) & = & C(s) \frac{\bar{r}_{ref}(s)}{\omega}\,,\quad
x_{ref}(0)=x_0,
\label{refu}\\
y_{ref}(t)&=&c^\top x_{ref}(t)\label{refy}\,,
\end{eqnarray}
where $\bar{r}_{ref}(s)$ is the Laplace transformation of the
signal
$$
\bar{r}_{ref}(t) = -\theta^{\top}(t) x_{ref}(t)-\sigma(t)+ k_g
r(t)\,,
$$
and $k_g$ is introduced in (\ref{kg_def}). The next  Lemma
establishes stability of the closed-loop system in
(\ref{refx})-(\ref{refy}).
\begin{lem}\label{lem:refstable}
If $D(s)$ verifies the condition in (\ref{condition3}), the
closed-loop reference system in (\ref{refx})-(\ref{refy})  is
stable.
\end{lem}
{\bf Proof.}
%
Let
\begin{equation}\label{Hsdef}
H(s)= (s I- A_m)^{-1} b\,.
\end{equation}
It follows from (\ref{refx})-(\ref{refy})  that
\begin{equation}\label{xrefs}
x_{ref}(s) = G(s) r_1(s) + H(s) C(s) k_g r(s)\,,
\end{equation}
where $r_1(s)$ is the Laplace transformation of
\begin{equation}\label{r1def}
r_1(t) = \theta^{\top}(t) x_{ref}(t)+\sigma(t)\,
\end{equation}
with the following bound:
\begin{equation}
\Vert r_1 \Vert_{\mathcal{L}_{\infty}} \leq L \Vert x_{ref}
\Vert_{\mathcal{L}_{\infty}}+\Vert \sigma
\Vert_{\mathcal{L}_{\infty}}\,.
\end{equation}
Since $D(s)$ verifies the condition in (\ref{condition3}), then
Theorem \ref{thm:sg}, applied to (\ref{xrefs}), ensures
 that the closed-loop system in (\ref{refx})-(\ref{refy})  is
 stable.
$\hfill{\square}$

\begin{lem}\label{lem:refstable2}
If $\theta(t)$ is  constant, and $\displaystyle{D(s)=1/s}$, then
the closed-loop reference system in (\ref{refx})-(\ref{refy}) is
stable {\it iff} the matrix $A_g$ in (\ref{condition4}) is
Hurwitz.
\end{lem}
{\bf Proof.}  In case of constant $\theta(t)$, the state space
form of the closed-loop system in (\ref{refx})-(\ref{refy}) is
given by:
\begin{equation}
\dot{x}_{ref}(t)  =  A_m x_{ref}(t)+b \left( \omega u_{ref}(t)+\theta^{\top} x_{ref}(t)+\sigma(t)\right), \label{xref_const}\\
\end{equation}
\begin{equation}
 \dot{u}_{ref}(t)  =  -\omega k u_{ref}(t) +k
\left(-\theta^{\top} x_{ref}(t)-\sigma(t)+ k_g r(t) \right),
\label{uref_const}
\end{equation}
\begin{equation}
 y_{ref}(t)=c^\top x_{ref}(t)\,.\label{yref_const}
\end{equation}
Letting
$$
\zeta(t) =\left[
\begin{array}{c}
x_{ref}(t)\\
 u_{ref}(t)
\end{array}\right]\,,
$$
it can be rewritten as
\begin{equation}\label{ref_5}
 \dot{\zeta}(t) = A_g \zeta(t) + \left[
\begin{array}{c}
b \sigma(t) \\
-k \sigma(t)+k k_g r(t)
\end{array}
\right]\,.
\end{equation}
We note that the LTI system in (\ref{ref_5}) is stable iff $A_g$
is Hurwitz, which concludes the proof. $\hfill{\square}$

\subsection{Bounded Error Signal}

\begin{lem}\label{lem:1}
For the system in (\ref{problemnew}) and the $\mathcal{L}_1$
adaptive controller in (\ref{L1_companionmodal}),
(\ref{adaptivelaw_L11})-(\ref{adaptivelaw_L13}) and
(\ref{controllaw}), the tracking error between the system state
and the companion model is bounded as follows:
\begin{equation}
\Vert \tilde{x} \Vert_{\mathcal{L}_{\infty}}  \leq
\sqrt{\frac{\theta_{m}}{\lambda_{\min}(P) \Gamma_c}}\,,
\label{barthemax}
\end{equation}
where
\begin{eqnarray}
 \quad \theta_{m}  \triangleq  \max_{\theta\in \Theta}
\sum_{i=1}^{n} 4 \theta_i^2+ 4 \Delta^2+4
\left(\omega_{u}-\omega_{l}
\right)^2& &\nonumber\\
\quad  +2\frac{\lambda_{\max}(P)}{\lambda_{\min}(Q)}\left(
  \max_{\theta\in \Theta}  \Vert\theta\Vert_2 d_{\theta} +   d_{\sigma}\Delta \right)
 \,.& &\label{thetamaxdef}
\end{eqnarray}
\end{lem}
{\bf Proof.}  Consider the following candidate Lyapunov function:
\begin{eqnarray}
V(\tilde{x}(t),\tilde{\theta}(t),\tilde{\omega}(t),\tilde{\sigma}(t))=
\tilde{x}^{\top}(t) P \tilde{x}(t)+ & & \nonumber\\
\Gamma_{c}^{-1} \tilde{\theta}^{\top}(t)
\tilde{\theta}(t)+\Gamma_{c}^{-1}
\tilde{\omega}^2(t)+\Gamma_{c}^{-1} \tilde{\sigma}^2(t)\,,& &
\nonumber
\end{eqnarray}
where
\begin{equation}\label{tildedef}
\tilde{\theta}(t) \triangleq \hat{\theta}(t)-\theta(t),\,
\tilde{\sigma}(t) \triangleq \hat{\sigma}(t)-\sigma(t),\,
\tilde{\omega}(t) \triangleq \hat{\omega}(t)-\omega\,.
\end{equation}
It follows from (\ref{problemnew}) and (\ref{L1_companionmodal})
that
\begin{equation}\label{error_dynamics}
\dot{\tilde{x}}(t)  =  A_m \tilde{x}(t)+b \left( \tilde{\omega}(t)
u(t) + \tilde{\theta}^{\top}(t) x(t)+\tilde{\sigma}(t) \right),\,
\tilde{x}(0)=0.
\end{equation}
Using the projection based adaptation laws from
(\ref{adaptivelaw_L11})-(\ref{adaptivelaw_L13}), one has the
following upper bound for $\dot V(t)$:
\begin{equation}\label{lem1_11}
\dot{V}(t) \leq -\tilde{x}^{\top}(t) Q \tilde{x}(t)
+\Gamma_{c}^{-1} \tilde{\theta}^{\top}(t)
\dot{\theta}(t)+\Gamma_{c}^{-1} \tilde{\sigma}
(t)\dot{\sigma}(t)\,.
\end{equation}
The projection algorithm ensures that $ \hat{\theta}(t) \in
\Theta$, $\hat{\omega}(t) \in \Omega$, $\hat{\sigma}(t) \in
\Delta$ for all $ t \geq 0, $ and therefore
\begin{eqnarray}
\max_{t\geq 0}\left( \Gamma_{c}^{-1} \tilde{\theta}^{\top}(t)
\tilde{\theta}(t)+\Gamma_{c}^{-1}
\tilde{\omega}^2(t)+\Gamma_{c}^{-1} \tilde{\sigma}^2(t)\right)\leq
& & \nonumber\\
 \left(\max_{\theta\in \Theta}
\sum_{i=1}^{n} 4 \theta_i^2+4 \Delta^2+4
\left(\omega_{u}-\omega_{l} \right)^2\right)/\Gamma_c & &
\label{lem1_5}
\end{eqnarray}
for any $t\geq 0$. If at any $t$
\begin{equation}\label{lem1_15}
V(t)>\frac{\theta_{m}}{\Gamma_c}\,,
\end{equation}
where $\theta_{m}$ is defined in (\ref{thetamaxdef}), then it
follows from (\ref{lem1_5}) that
\begin{equation}
\tilde{x}^{\top}(t) P \tilde{x}(t) >2
\frac{\lambda_{\max}(P)}{\Gamma_c \lambda_{\min}(Q)}\left(
  \max_{\theta\in \Theta}  \Vert\theta\Vert_2 d_{\theta} +
  d_{\sigma}\Delta
\right),
\end{equation}
and hence
\begin{eqnarray}
\tilde{x}^{\top}(t) Q \tilde{x}(t)  >  \frac{\lambda_{\min}(Q)}{
\lambda_{\max}(P)} \tilde{x}^{\top}(t) P
\tilde{x}(t) & & \nonumber\\
 > 2 \frac{
  \max_{\theta\in \Theta}  \Vert\theta\Vert_2 d_{\theta} +   d_{\sigma}\Delta
 }{\Gamma_c}\,.& & \nonumber
\end{eqnarray}
The upper bounds in  (\ref{derivativeB}) along with the projection
based adaptive laws lead to the following upper bound:
\begin{equation}
 \frac{ \tilde{\theta}^{\top}(t)
\dot{\theta}(t)+ \tilde{\sigma}(t) \dot{\sigma}(t)}{\Gamma_c} \leq
2\frac{\max_{\theta\in \Theta}  \Vert\theta\Vert_2 d_{\theta} +
 d_{\sigma}\Delta
  }{\Gamma_c}\,.
\end{equation}
Hence,  if $\displaystyle{V(t)>\frac{\theta_{m}}{\Gamma_c}}$, then
from (\ref{lem1_11})  we have
\begin{equation}\label{lem1_21}
\dot{V}(t) <0\,.
\end{equation}
 Since we have set $\hat{x}(0)=x(0)$, we can verify that
\begin{eqnarray}
V(0) \leq  \Big(\max_{\theta\in \Theta} \sum_{i=1}^{n} 4
\theta_i^2+ 4 \Delta^2+4 \left(\omega_{u}-\omega_{l}
\right)^2\Big)/\Gamma_c
  <
\frac{\theta_{m}}{\Gamma_c}\,.\qquad\qquad\qquad\qquad\qquad\qquad\qquad\qquad\qquad
& & \nonumber
\end{eqnarray}
It follows from (\ref{lem1_21}) that $\displaystyle{ V(t) \leq
\frac{\theta_{m}}{\Gamma_c}}$ for any $t \geq 0$. Since $
\lambda_{\min}(P) \Vert \tilde{x}(t) \Vert^2 \leq
\tilde{x}^{\top}(t) P \tilde{x}(t)\leq V(t)$, then
$$ 
\displaystyle{|| \tilde{x}(t) ||^2  \leq
\frac{\theta_{m}}{\lambda_{\min}(P) \Gamma_c}}\,,
$$ 
which concludes the proof.
 $\hfill{\square}$

\begin{rem}
We note that the bound in (\ref{barthemax}) is similar to the
bounds derived in \cite{kkk}, assuming appropriate trajectory
initialization  to ensure transient performance improvement for
system's output tracking. For the particular control architecture
in this paper, the appropriate trajectory initialization is
ensured by setting
 $\hat{x}(0)=x(0)$. However, due to the special filtering technique
subject to $\mathcal L_1$-gain requirement, we obtain uniform
smooth transient for systems's both signals, input and output, as
proved in the next section.
\end{rem}


\subsection{Transient Performance}
Let
\begin{equation}\label{HsDefnew}
H(s) = (s I -A_m)^{-1} b\,.
\end{equation}
It follows from Lemma \ref{lem:pre2} that there exists $c_{o}\in
{\rr}^n$ such that
\begin{equation}\label{thm5_pre}
c_{o}^{\top} H(s) = \frac{N_n(s)}{N_d(s)}\,,
\end{equation}
where the order of $N_d(s)$ is one more than the order of
$N_n(s)$, and both $N_n(s)$ and $N_d(s)$ are stable polynomials.

\begin{thm}\label{thm:5}
Given the system in (\ref{problemnew}) and the $\mathcal{L}_1$
adaptive controller defined via (\ref{L1_companionmodal}),
(\ref{adaptivelaw_L11})-(\ref{adaptivelaw_L13}) and
(\ref{controllaw}) subject to (\ref{condition3}), we have:
\begin{eqnarray}
\Vert x-x_{ref} \Vert_{{\mathcal L}_{\infty}}  & \leq & \gamma_1\,
,\label{thm5_00} \\
\Vert u - u_{ref} \Vert_{{\mathcal L}_{\infty}} & \leq &
\gamma_2\, ,\label{thm5_01}
\end{eqnarray}
where
\begin{eqnarray}
\gamma_1 & = & \frac{\Vert C(s) \Vert_{\mathcal{L}_1}}{1-\Vert
H(s) (1-C(s))\Vert_{\mathcal{L}_1} L}
\sqrt{\frac{\theta_{m}}{\lambda_{\max}(P) \Gamma_c}}\,,
\label{gamma1def} \\
\gamma_2 & = & \left\Vert \frac{C(s)}{\omega}
\right\Vert_{\mathcal{L}_1} L \gamma_1 + \nonumber\\
& & \Big\| \frac{C(s)}{\omega} \frac{1}{c_{o}^{\top} H(s)}
c_{o}^{\top} \Big\|_{\mathcal{L}_1}
\sqrt{\frac{\theta_{m}}{\lambda_{\max}(P) \Gamma_c}}\,.
\label{gamma2def}
\end{eqnarray}
\end{thm}
\medskip

{\bf Proof.}
Let
\begin{eqnarray}
\tilde{r}(t) & = & \tilde{\omega}(t) u(t)
+\tilde{\theta}^{\top}(t) x(t) + \tilde{\sigma}(t)\,,\nonumber\\
r_2(t) & = & \theta^{\top}(t) x(t)+\sigma(t)\,.\nonumber
\end{eqnarray}
It follows from (\ref{controllaw}) that
$$
\chi(s) =D(s) (\omega u(s) + r_2(s)-k_g r(s)+\tilde{r}(s))\,,
$$
where $\tilde r(s)$ and $r_2(s)$ are the Laplace transformations
of signals $\tilde r(t)$ and $r_2(t)$. Consequently
\begin{eqnarray}
\chi(s) & = & \frac{D(s)}{1+k\omega  D(s)}
(r_2(s)-k_g r(s)+\tilde{r}(s))\,,\label{thm1_21}\\
 u(s) & = & -\frac{k D(s)}{1+k
\omega D(s)} (r_2(s)-k_g r(s)+\tilde{r}(s))\,.\label{thm1_22}
\end{eqnarray}
Using the definition of $C(s)$ from  (\ref{Csdef}), we can write
\begin{equation}\label{thm1_22new}
\omega u(s) = -C(s) (r_2(s)-k_g r(s)+\tilde{r}(s))\,,
\end{equation}
and the system in (\ref{problemnew}) consequently takes the form:
\begin{equation}\label{thm1_5}
x(s) = H(s)  \left( (1-C(s)) r_2(s) +C(s) k_g r(s)-  C(s)
\tilde{r}(s)\right).
\end{equation}
It follows from (\ref{refx})-(\ref{refu}) that
\begin{equation}\label{thm1_6}
x_{ref}(s) = H(s) \left( (1-C(s)) r_1(s) +C(s) k_g r(s)\right)\,,
\end{equation}
where $r_1(s)$ is the Laplace transformation of the signal
$r_1(t)$ defined in (\ref{r1def}). Let $e(t)=x(t)-x_{ref}(t)$.
Then, using (\ref{thm1_5}), (\ref{thm1_6}), one gets
\begin{equation}
e(s) = H(s) \left( (1-C(s)) r_3(s) -   C(s) \tilde{r}(s)\right),
e(0)=0\,,\label{thm1_7}
\end{equation}
where $r_3(s)$ is the Laplace transformation of the signal
\begin{equation}\label{thm1_11}
r_3(t) = \theta^{\top}(t) e(t)\,.
\end{equation}
Lemma \ref{lem:1} gives the following upper bound:
\begin{equation}\label{thm1_8}
\Vert e_{t} \Vert_{\mathcal{L}_{\infty}} \leq \Vert H(s)
(1-C(s))\Vert_{\mathcal{L}_1}  \Vert r_{3_t}
\Vert_{\mathcal{L}_{\infty}} + \Vert r_{4_t}
\Vert_{\mathcal{L}_{\infty}}\,,
\end{equation}
where $ r_4(t)$ is the signal with its Laplace transformation
$$
r_4(s) = C(s) H(s) \tilde{r}(s).
$$
From the relationship in  (\ref{error_dynamics}) we have
\begin{equation}
\tilde{x}(s) = H(s) \tilde{r}(s)\,,
\end{equation}
which leads to
\begin{equation}
r_4(s) = C(s) \tilde{x}(s)\,,
\end{equation}
and hence
\begin{equation}\label{thm1_9}
\Vert r_{4_t} \Vert_{\mathcal{L}_{\infty}} \leq \Vert C(s)
\Vert_{\mathcal{L}_1} \Vert \tilde{x}_t
\Vert_{\mathcal{L}_{\infty}}\,.
\end{equation}
Using the definition of $L$ in (\ref{thetabardef}), one can verify
 easily that
\begin{equation}\label{Lproperty}
\Vert (\theta^{\top} e)_t \Vert_{\mathcal{L}_{\infty}} \leq L
\Vert
 e_t \Vert_{\mathcal{L}_{\infty}}\,,
\end{equation}
and hence  the following upper bound can be derived from
(\ref{thm1_11}):
\begin{equation}\label{thm1_12}
\Vert r_{3_t} \Vert_{\mathcal{L}_{\infty}} \leq L \Vert e_{t}
\Vert_{\mathcal{L}_{\infty}}\,.
\end{equation}
From (\ref{thm1_8}) we have
\begin{equation}
\Vert e_{t} \Vert_{\mathcal{L}_{\infty}} \leq \Vert H(s)
(1-C(s))\Vert_{\mathcal{L}_1} L \Vert e_{t}
\Vert_{\mathcal{L}_{\infty}}+  \Vert C(s)
\Vert_{\mathcal{L}_1}\Vert \tilde{x}_t
\Vert_{\mathcal{L}_{\infty}}\,.
\end{equation}
The upper bound from Lemma \ref{lem:1} and the $\mathcal L_1$-gain
requirement from (\ref{condition3}) lead to the following upper
bound
\begin{equation}
\Vert e_{t} \Vert_{\mathcal{L}_{\infty}} \leq \frac{\Vert C(s)
\Vert_{\mathcal{L}_1}}{1-\Vert H(s) (1-C(s))\Vert_{\mathcal{L}_1}
L } \sqrt{\frac{\theta_{m}}{\lambda_{\max}(P) \Gamma_c}}\,,
\end{equation}
which holds uniformly for all $t\ge 0$ and therefore leads to
(\ref{thm5_00}).

To prove the bound in (\ref{thm5_01}), we notice that from
(\ref{refu}) and (\ref{thm1_22new})  one can derive
\begin{equation}\label{lem2_new5}
u(s)-u_{ref}(s) = - \frac{C(s)}{\omega} \theta^{\top}(t)
(x(s)-x_{ref}(s)) -r_5(s)\,,
\end{equation}
where $r_5(s)= \frac{C(s)}{\omega}\tilde{r}(s)$. Therefore, it
follows from Lemma \ref{lem:1} that
\begin{equation}\label{thm1_99}
\Vert u-u_{ref}\Vert_{\mathcal{L}_{\infty}} \leq \frac{\Vert C(s)
\Vert_{\mathcal{L}_1} L}{\omega} \Vert
x-x_{ref}\Vert_{\mathcal{L}_{\infty}}+\Vert
r_5\Vert_{\mathcal{L}_{\infty}}\,.
\end{equation}
We have
\begin{eqnarray}
& & r_5(s)  = \frac{C(s)}{\omega} \frac{1}{c_{o}^{\top} H(s)}
c_{o}^{\top} H(s)
\tilde r(s)      \nonumber\\
& & \qquad =  \frac{C(s)}{\omega} \frac{1}{c_{o}^{\top} H(s)}
c_{o}^{\top} \tilde{x}(s)\,,\nonumber
\end{eqnarray}
where $c_{o}$ is introduced in (\ref{thm5_pre}). Using the
polynomials from (\ref{thm5_pre}), we can write that
$$ 
 \frac{C(s)}{\omega} \frac{1}{c_{o}^{\top} H(s)} = \frac{C(s)}{\omega}
 \frac{N_d(s)}{N_n(s)}\,,
$$ 
where $N_d(s)$, $N_n(s)$ are stable polynomials and the order of
$N_n(s)$ is one less than the order of $N_d(s)$.  Since $C(s)$ is
stable and strictly proper, the complete system $
 C(s) \frac{1}{c_{o}^{\top} H(s)}
$ is proper and stable, which implies that its $\mathcal{L}_1$
gain exists and is finite. Hence, we have
$$ 
\Vert r_5\Vert_{{\mathcal L}_{\infty}} \leq \Big\|
\frac{C(s)}{\omega}  \frac{1}{c_{o}^{\top} H(s)} c_{o}^{\top}
\Big\|_{\mathcal{L}_1} \Vert \tilde{x} \Vert_{{\mathcal
L}_{\infty}}\,.
$$ 
 Lemma \ref{lem:1}  consequently leads to the upper bound:
$$ 
\Vert r_5\Vert_{{\mathcal L}_{\infty}} \leq  \Big\|
\frac{C(s)}{\omega} \frac{1}{c_{o}^{\top} H(s)} c_{o}^{\top}
\Big\|_{\mathcal{L}_1} \sqrt{\frac{\theta_{m}}{\lambda_{\max}(P)
\Gamma_c}}\,,
$$ 
which, when substituted into (\ref{thm1_99}), leads to
(\ref{thm5_01}).
 $\hfill{\square}$

\begin{thm}\label{thm:2}
 For the closed-loop system in
(\ref{problemnew}) with $\mathcal{L}_1$ adaptive controller
defined via (\ref{L1_companionmodal}),
(\ref{adaptivelaw_L11})-(\ref{adaptivelaw_L13}) and
(\ref{controllaw}), subject to (\ref{condition4}), if $\theta(t)$
is (unknown) constant and $\displaystyle{D(s)=\frac{1}{s}}$, we
have:
\begin{eqnarray}
\Vert x-x_{ref} \Vert_{{\mathcal L}_{\infty}}  & \leq & \gamma_3\,
,\label{thm2_00} \\
\Vert u - u_{ref} \Vert_{{\mathcal L}_{\infty}} & \leq &
\gamma_4\, ,\label{thm2_11}
\end{eqnarray}
where
\begin{eqnarray}
\gamma_3 & = & \Big\| H_g(s) C(s) \frac{1}{c_{o}^{\top} H(s)}
c_{o}^{\top} \Big\|_{\mathcal{L}_1}
\sqrt{\frac{\theta_{m}}{\lambda_{\max}(P) \Gamma_c}},
\label{gamma3def} \\
\gamma_4 & = & \left\Vert \frac{C(s)}{\omega} \theta^{\top}
\right\Vert_{\mathcal{L}_1} \gamma_3 + \nonumber\\
& & \Big\| \frac{C(s)}{\omega} \frac{1}{c_{o}^{\top} H(s)}
c_{o}^{\top} \Big\|_{\mathcal{L}_1}
\sqrt{\frac{\theta_{m}}{\lambda_{\max}(P) \Gamma_c}}\,,
\label{gamma4def}
\end{eqnarray}
and $$ H_g(s)=(s I - A_g) \left[
\begin{array}{c}
b \\
0
\end{array}
\right].
$$
\end{thm}
{\bf Proof.} Recall that for constant $\theta$ we had

$$D(s)=\frac{1}{s},\quad C(s)=\frac{k\omega}{s+k \omega}.$$
Let
$$
\zeta(s) = -\frac{C(s)}{\omega} \theta^{\top} e(s)\,.
$$
With this notation, (\ref{thm1_7}) can be written as
$$
e(s)=H(s) \left( \theta^{\top} e(s) +\omega \zeta(s) -C(s)
\tilde{r}(s) \right)
$$
and further put into state space form as:
\begin{equation}
\left[
\begin{array}{c}
\dot{e}(t) \\
\dot{\zeta}(t)
\end{array}
\right]= A_g \left[\begin{array}{c}
e(t) \\
\zeta(t)
\end{array}
\right] + \left[
\begin{array}{c}
b \\
0
\end{array}
\right] r_6(t)\,,\label{thm2_5}
\end{equation}
where $r_6(t)$ is the signal with its Laplace transformation
\begin{equation}
r_6(s) = -C(s) \tilde{r}(s)\,.
\end{equation}
Let
$$
x_{\zeta}(t) =[e^{\top}(t)\,\, \zeta(t)]^{\top}.
$$
 Since $A_g$ is Hurwitz, then $H_g(s)$ is stable and strictly
proper. It follows from (\ref{thm2_5}) that
$$
x_{\zeta}(s) = -H_g(s) C(s)  \tilde{r}(s)\,.
$$
Therefore, we have
\begin{eqnarray}
& &x_{\zeta}(s)  = -H_g(s) C(s) \frac{1}{c_{o}^{\top} H(s)}
c_{o}^{\top} H(s)
\tilde r(s)      \nonumber\\
& & \qquad =  -H_g(s) C(s) \frac{1}{c_{o}^{\top} H(s)}
c_{o}^{\top} \tilde{x}(s)\,,\nonumber
\end{eqnarray}
where $c_{o}$ is introduced in (\ref{thm5_pre}).
 It follows from (\ref{thm5_pre}) that
$ 
 H_g(s) C(s) \frac{1}{c_{o}^{\top} H(s)} = H_g(s) C(s)
 \frac{N_d(s)}{N_n(s)}\,,
$ 
where $N_d(s)$, $N_n(s)$ are stable polynomials and the order of
$N_n(s)$ is one less than the order of $N_d(s)$.  Since both
$H_g(s)$ and $C(s)$ are stable and strictly proper, the complete
system $
 H_g(s) C(s) \frac{1}{c_{o}^{\top} H(s)}
$ is proper and stable, which implies that its $\mathcal{L}_1$
gain exists and is finite. Hence, we have
$$ 
\Vert x_{\zeta}\Vert_{{\mathcal L}_{\infty}} \leq \Big\| H_g(s)
C(s) \frac{1}{c_{o}^{\top} H(s)} c_{o}^{\top}
\Big\|_{\mathcal{L}_1} \Vert \tilde{x} \Vert_{{\mathcal
L}_{\infty}}\,.
$$ 

The proof of (\ref{thm2_11}) is similar to the  proof of
(\ref{thm5_01}). $\hfill{\square}$

\begin{cor}\label{cor:1}
Given the system in (\ref{problemnew}) and the $\mathcal{L}_1$
adaptive controller defined via (\ref{L1_companionmodal}),
(\ref{adaptivelaw_L11})-(\ref{adaptivelaw_L13}) and
(\ref{controllaw}) subject to (\ref{condition3}), we have:
\begin{eqnarray}
\lim_{\Gamma_c\rightarrow \infty} \left(x(t)-x_{ref}(t) \right)& =
& 0\, ,\qquad \forall t\geq 0,\label{cor1_00}
\\
\lim_{\Gamma_c\rightarrow \infty} \left(u(t)-u_{ref}(t)\right) & =
& 0\,,\qquad \forall t\geq 0\,.\label{cor1_01}
\end{eqnarray}
\end{cor}

Thus, the tracking error between $x(t)$ and $x_{ref}(t)$, as well
between $u(t)$ and  $u_{ref}(t)$, is uniformly bounded by a
constant inverse proportional to $\Gamma_c$. This implies that
during the transient one can achieve arbitrarily close tracking
performance for both signals simultaneously by   increasing
$\Gamma_c$.


\subsection{Asymptotic Convergence}
Since the bounds in (\ref{thm5_00}) and (\ref{thm5_01}) are
uniform for all $t\ge 0$, they are in charge for both transient
and steady state performance.  In case of constant $\theta$ one
can prove in addition the following asymptotic result.

\begin{lem}\label{lem:5}
Given the system in (\ref{problemnew}) with constant $\theta$  and
$\mathcal{L}_1$ adaptive controller defined via
(\ref{L1_companionmodal}),
(\ref{adaptivelaw_L11})-(\ref{adaptivelaw_L13}) and
(\ref{controllaw}) subject to (\ref{condition3}), we have:
\begin{eqnarray}
\lim_{t\rightarrow \infty} \tilde{x}(t) & = & 0\,.
\end{eqnarray}
\end{lem}
{\bf Proof:} It follows from Lemmas \ref{lem:refstable} and
\ref{lem:1}, and Theorem \ref{thm:5} that both $x(t)$ and
$\hat{x}(t)$ in $\mathcal{L}_1$ adaptive controller are bounded
for bounded reference inputs. The adaptive laws in
(\ref{adaptivelaw_L11})-(\ref{adaptivelaw_L13})  ensure that the
estimates $\hat{\theta}(t)$, $\hat{\omega}(t)$, $\hat{\sigma}(t)$
are also bounded. Hence, it can be checked easily from
(\ref{error_dynamics}) that $\dot{\tilde{x}}(t)$ is bounded, and
it follows from Barbalat's lemma that $\displaystyle{
\lim_{t\rightarrow \infty} \tilde{x}(t) =0} $. $\hfill{\square}$

%
%
%

\subsection{Design Guidelines}

We note that the control law $u_{ref}(t)$ in the closed-loop
reference system, which is used in the analysis of $\mathcal
L_\infty$ norm bounds, is not implementable since its definition
involves the unknown parameters.  Theorem \ref{thm:5} ensures that
the $\mathcal L_1$ adaptive controller approximates  $u_{ref}(t)$
both in transient and steady state. So, it is important to
understand how these bounds can be used for ensuring uniform
transient response with {\em desired} specifications. We notice
that the following {\em ideal} control signal
\begin{equation}\label{uideal}
u_{ideal}(t)=\frac{k_g r(t)-\theta^{\top}(t) x_{ref}(t)-\sigma(t)
}{\omega}
\end{equation}
is the one that leads  to desired system response:
\begin{eqnarray}
\dot{x}_{ref}(t) &=& A_m x_{ref}(t)+b k_g r(t)\label{desrefstate}\\
y_{ref}(t)&=&c^\top x_{ref}(t)\label{desrefoutput}
\end{eqnarray}
by cancelling the uncertainties exactly. In the closed-loop
reference system (\ref{refx})-(\ref{refy}),
 $u_{ideal}(t)$ is further low-pass filtered by $C(s)$  in
(\ref{refu}) to have guaranteed low-frequency range. Thus, the
reference system in (\ref{refx})-(\ref{refy}) has a different
response as compared to (\ref{desrefstate}), (\ref{desrefoutput})
with (\ref{uideal}). In \cite{CaoHovakimyan1}, specific design
guidelines are suggested for selection of $C(s)$ to ensure that in
case of constant $\theta$ the response of (\ref{xref_const}),
(\ref{uref_const}), (\ref{yref_const}) can be made as close as
possible to (\ref{desrefstate}), (\ref{desrefoutput}) with
(\ref{uideal}). In case of fast varying $\theta(t)$, it is obvious
that the bandwidth of the controller needs to be matched
correspondingly.

\section{Simulations}\label{sec:simu}
As an illustrative example, consider a single-link robot arm which
is rotating on a vertical plane. The system dynamics are given by:
\begin{equation}\label{simu_problem}
I \ddot{q}(t) + \frac{M g L \cos q(t)}{2}+F(t) \dot{q}(t)+F_1(t)
q(t)+\bar{\sigma}(t) =u(t)\,,
\end{equation}
where $q(t)$ and $\dot{q}(t)$ are measured angular position and
velocity, respectively, $u(t)$ is the input torque, $I$ is the
unknown moment of inertia, $M$ is the unknown mass, , $L$ is the
unknown length, $F(t)$ is an unknown time-varying friction
coefficient, $F_1(t)$ is position dependent external torque, and
$\bar{\sigma}(t)$ is unknown bounded disturbance. The control
objective is to design $u(t)$ to achieve tracking of bounded
reference input $r(t)$ by $q(t)$. Let
$$
x = [q\quad \dot{q}]^{\top}\,.
$$
The system in (\ref{simu_problem}) can be presented in the
state-space form as:
\begin{eqnarray}
\dot{x}(t) & = & A x(t) + b \Big( \frac{u(t)}{I}  +\frac{M g L
\cos(x_1(t))}{2 I}+\frac{\bar{\sigma}(t)}{I}\nonumber\\
& & +\frac{F_1(t)}{I} x_1(t)+\frac{F(t)}{I} x_2(t)\Big),\quad x(0)=x_0 \,,\nonumber\\
y(t) & = & c^{\top} x(t)\,,\label{simu_pro2}
\end{eqnarray}
where $x_0$ is the initial condition,
\begin{equation}
A= \left[
\begin{array}{cc}
0 &  1\\
0 & 0
\end{array}
\right]\,,\qquad b=\left[
\begin{array}{c}
0 \\
1
\end{array}
\right]\,, \qquad c=\left[
\begin{array}{c}
1 \\
0
\end{array}
\right]\,.\label{Abdef}
\end{equation}
\begin{figure}[h!]
\centering
\mbox{{\epsfig{file=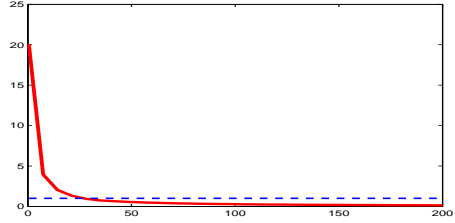,width=2.8in,height=1.3in
}\label{fig:criterion}}} \caption{$\Vert G(s)
\Vert_{\mathcal{L}_1} L$
 with respect to $\omega k$.}
\end{figure} %
%

The system can be further put into the form:
\begin{eqnarray}
\dot{x}(t) & = & A_m x(t) + b (\omega u(t)+\theta^{\top}(t) x(t)+ \sigma(t))\,,\nonumber\\
y(t) & = & c^{\top} x(t)\,,\quad x(0)=x_0 \,,\nonumber
\end{eqnarray}
where $\omega  = \frac{1}{I}$ is the unknown control
effectiveness,
\begin{eqnarray}
A_m & = & \left[
\begin{array}{lr}
0 & 1 \\
-1 & -1.4
\end{array}
\right], b  =  \left[
\begin{array}{l}
0\\
1
\end{array}
\right]\,, c  = \left[
\begin{array}{l}
1\\
0
\end{array}
\right]\,,
\end{eqnarray}
\begin{eqnarray}
\theta(t) & = & \left[1+\frac{F_1(t)}{I}\,\;\;\; \, 1.4+\frac{F(t)}{I}\right]^{\top}\,,\nonumber\\
\sigma(t) & = & \frac{M g L \cos(x_1(t))}{2 I}
 +\frac{\bar{\sigma}(t)}{I}\,.\nonumber
\end{eqnarray}
Let the unknown control effectiveness, time-varying parameters and
disturbance be given by:
\begin{eqnarray}
\omega & = & 1\,,\nonumber\\
\theta(t) & = & [2+ \cos(\pi t)\;\;\; 2+0.3\sin(\pi t)+0.2\cos(2
t)]^{\top}\,,\nonumber\\
\sigma(t) & = & \sin(\pi t)\,,
\end{eqnarray}
so that the compact sets can be conservatively chosen as
\begin{equation}\label{compsets}
\Omega=[0.2,\,5],\, \Theta=[-10,\,10],\, \Delta=[-10,\,10]\,.
\end{equation}
\begin{figure}[h!]
\centering \mbox{\subfigure[ $x_1(t)$ (solid),  $\hat{x}_1(t)$
(dashed), and $r(t)$(dotted)
 ]{\epsfig{file=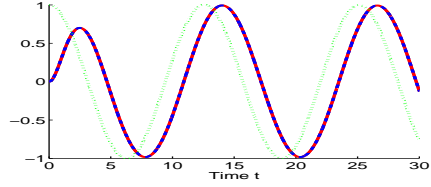,width=2.5in,height=1.in
}\label{fig:L1_y}}} \mbox{ \subfigure[ Time-history of $u(t)$
]{\epsfig{file=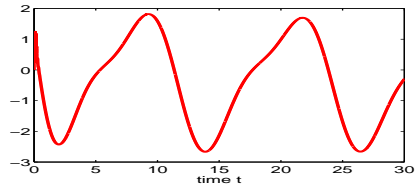,width=2.5in, height=1.in}\label{fig:L1_u}
}} \caption{Performance of $\mathcal{L}_1$ adaptive controller for
$\sigma(t)=\sin(\pi t)$}
\end{figure} %
\begin{figure}[h!]
\centering \mbox{\subfigure[ $x_1(t)$ (solid),  $\hat{x}_1(t)$
(dashed), and $r(t)$(dotted)
 ]{\epsfig{file=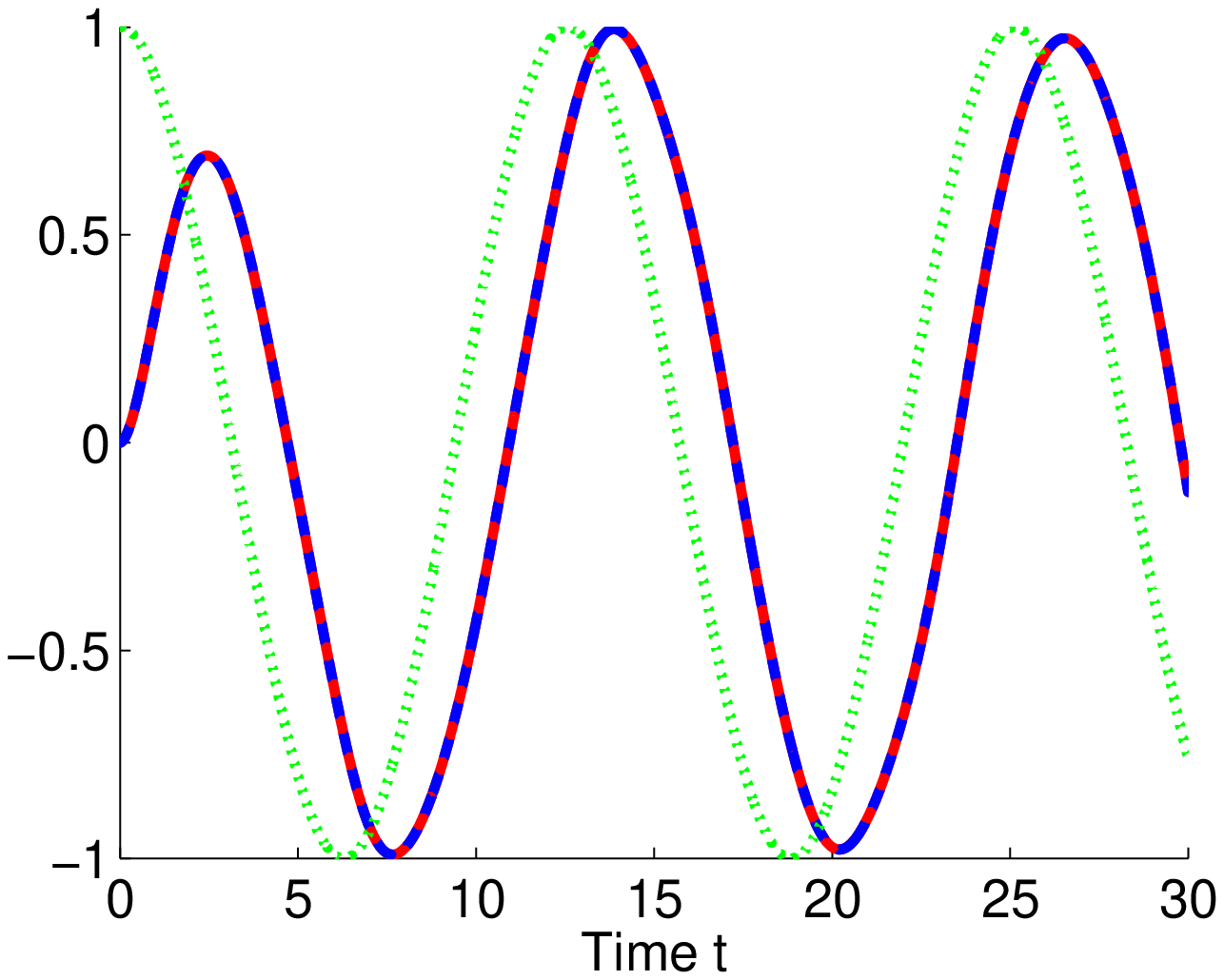,width=2.5in,height=1.in
}\label{fig:L1_yd}}} \mbox{ \subfigure[ Time-history of $u(t)$
]{\epsfig{file=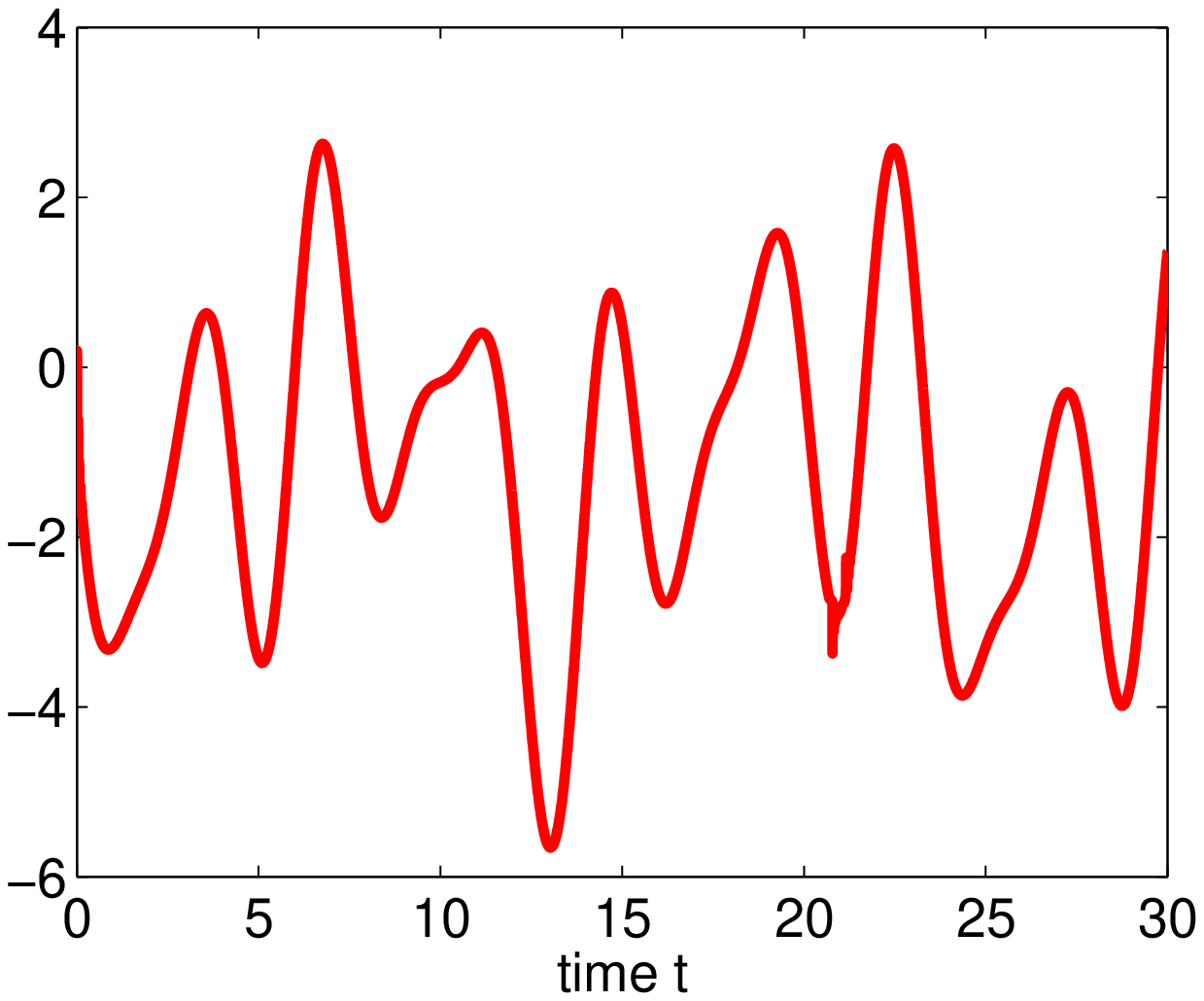,width=2.5in,
height=1.in}\label{fig:L1_ud} }} \caption{Performance of
$\mathcal{L}_1$ adaptive controller for $\sigma(t)=\cos(x_1(t))+ 2
\sin(10 t)+\cos(15  t)$}
\end{figure} %
\begin{figure}[h!]
\centering \mbox{\subfigure[ $x_1(t)$ (solid),  $\hat{x}_1(t)$
(dashed), and $r(t)$(dotted)
 ]{\epsfig{file=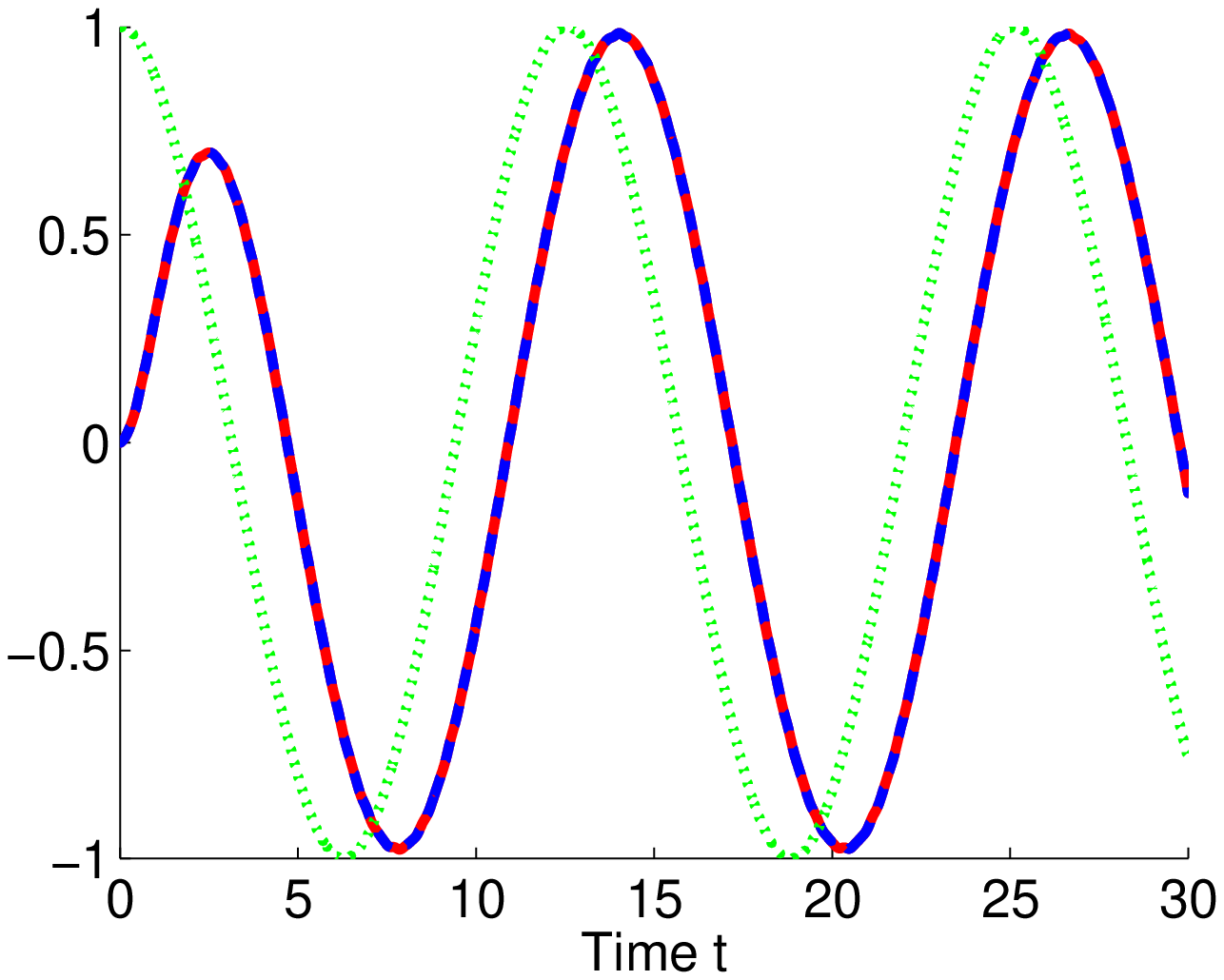,width=2.5in,height=1.in
}\label{fig:L1_ydh}}} \mbox{ \subfigure[ Time-history of $u(t)$
]{\epsfig{file=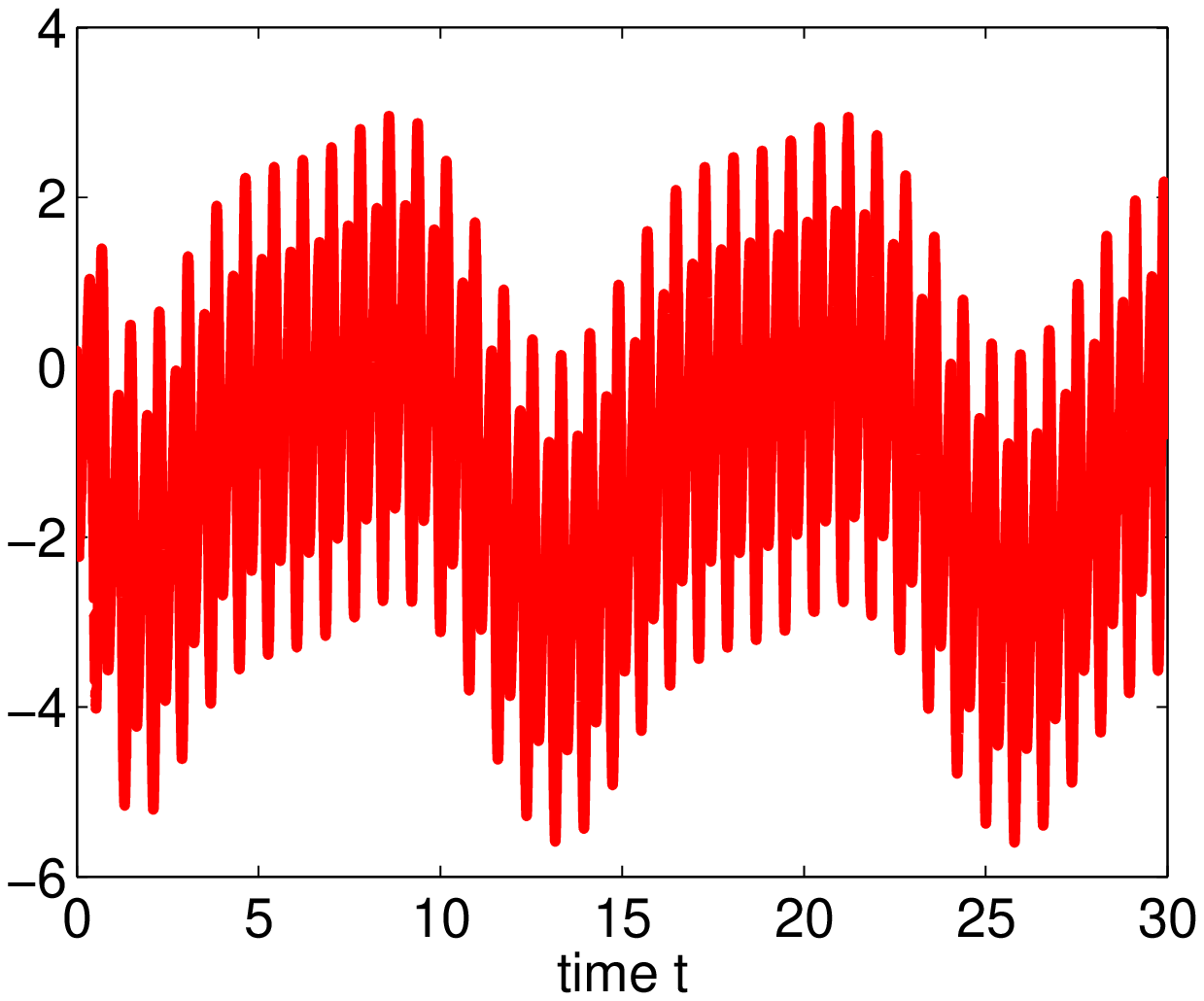,width=2.5in,
height=1.in}\label{fig:L1_udh} }} \caption{Performance of
$\mathcal{L}_1$ adaptive controller for $\sigma(t)=\cos(x_1(t))+ 2
\sin(100  t)+\cos(150 t)$}
\end{figure} %
For implementation of the $\mathcal{L}_1$ adaptive controller
(\ref{L1_companionmodal}),
(\ref{adaptivelaw_L11})-(\ref{adaptivelaw_L13}) and
(\ref{controllaw}), we need to verify the $\mathcal{L}_1$
stability requirement in (\ref{condition3}). Letting
$$
D(s)=1/s\,,
$$
we have
\begin{equation}\label{simu_ydes}
G(s) = \frac{\omega k}{s+\omega k}  H(s),
\end{equation}
where
\begin{equation}
H(s) = \left[
\begin{array}{c}
\frac{1}{s^2+1.4 s+1} \\
\frac{s}{s^2+1.4 s+1}
\end{array}
\right]\,.
\end{equation}
We can check easily that for our selection of compact sets in
(\ref{compsets}), the resulting $L=20$ in (\ref{thetabardef}). In
Fig. \ref{fig:criterion}, we plot $\Vert G(s)
\Vert_{\mathcal{L}_1} L$ as a function of $\omega k$ and compare
it to $1$. We notice that for $\omega k
> 30$, we have $\Vert G(s) \Vert_{\mathcal{L}_1} L<1$.
Since $\omega>0.5$, we set $k=60$. At last, we set the adaptive
gain as $\Gamma_c = 10000 $.

The simulation results of the $\mathcal{L}_1$ adaptive controller
are shown in Figures \ref{fig:L1_y}-\ref{fig:L1_u} for reference
input $ r=\cos(\pi t)$. Next, we consider different disturbance
signal:
$$
\sigma(t)  =  \cos(x_1(t))+ 2 \sin(10 t)+\cos(15  t)\,.
$$
The simulation results are shown in
\ref{fig:L1_yd}-\ref{fig:L1_ud}. Finally, we consider much higher
frequencies in the disturbance:
$$
\sigma(t)  =  \cos(x_1(t))+ 2 \sin(100  t)+\cos(150 t)\,.
$$
The simulation results are shown in
\ref{fig:L1_ydh}-\ref{fig:L1_udh}. We note that the
$\mathcal{L}_1$ adaptive controller guarantees smooth and uniform
transient performance in the presence of different unknown
nonlinearities and time-varying disturbances. The controller
frequencies are exactly matched with the frequencies of the
disturbance that it is supposed to cancel out. We also notice that
$x_1(t)$ and $\hat{x}_1(t)$ are almost the same in Figs.
\ref{fig:L1_y}, \ref{fig:L1_yd} and \ref{fig:L1_ydh}.

\section{Conclusion}\label{sec:con}
A novel $\mathcal{L}_1$ adaptive control architecture is presented
 that has guaranteed transient response in addition
to stable tracking for systems with time-varying unknown
parameters and bounded disturbances. The control signal and the
system response approximate the same signals of a closed-loop
reference LTI system, which can be designed to achieve desired
specifications. In Part II of this paper
\cite{CDC06_chengyu_big_L1}, we derive the stability margins of
this $\mathcal{L}_1$ adaptive control architecture.


\end{document}